\input amstex
\documentstyle {amsppt}
\pageheight{50.5pc}
\pagewidth{32pc}
\define\E{{\bold E}}
\redefine\P{{\bold P}}
\define\sign{\,\hbox{{\rm sign}}\,}
\topmatter
\title{On a multidimentional Brownian motion with partly reflecting
 membrane on a hyperplane}
\endtitle
\author
 Ludmila L. Zaitseva
\endauthor
\address
 Institute of Mathematics, Kiev, Ukraine
\endaddress
\abstract
A multidimensional Brownian motion with partial reflection on a hyperplane $S$
in the direction $qN+\alpha $, where $N$ is the conormal vector to the
hyperplane and $q\in [-1,1], \alpha \in S$ are given parametres, is constructed
and this construction is based on both analytic and probabilistic approaches. The
joint distribution of $d-$dimensional analogy to skew Brownian motion and its local
time on the hyperplane is obtained.
\endabstract
\email
zaitseva\@univ.kiev.ua
\endemail
\endtopmatter
\document
\rightheadtext{Brownian motion with partly reflecting membrane}
\leftheadtext{Zaitseva L.L.}

\head
{\S 1. Introduction.}
\endhead

In this paper we construct a generalized diffusion process in a
 $d$-dimensional Euclidean space $\Re ^d$ for which the diffusion matrix $B$
 is constant and the drift vector is equal to $(qB\nu +\alpha )\delta _S(x)$
 where $\nu \in \Re ^d$ is a given unit vector, $S$ is the hyperplane in $\Re
 ^d$ orthogonal to $\nu $, the parameters $q\in [ -1,1]$ and $\alpha \in S$
 are given  and $\delta _S(x)$ is a generalized function that is determined
 by the relation
$$
\int_{\Re ^d}{\varphi (y)\delta _S(y)dy}=\int_{S}{\varphi (y) d\sigma }
$$
valid for any test function $\varphi $.

This problem in a more general case (namely, the parameters $\alpha ,q$ are
some functions of $x\in S$ and the diffusion matrix $B$ is a function of
$x\in \Re ^d$) was considered by Kopytko B.I. \cite{1} and solved by an
analytic method. More simple case ($B$ is an identical operator, $\alpha ,q$
are some constants) was considered by Kopytko B.I., Portenko N.I. \cite{2}
with use of an analytic method too.

In this work a probabilistic approach to this problem will be proposed. We
construct the process desired as a solution of the following stochastic
differential equation
$$
dx(t)=(qB\nu +\alpha )\delta _S(x(t))dt + B^{1/2}dw(t)
$$
where $B^{1/2}$ is the positive square root from the operator $B$ and $w(t)$
is a given Wiener process in $\Re ^d$. If $\alpha =0$ then the solution of
this equation is known. Let us denote it by $\tilde{x}(t)$. The idea is to
add the process $\alpha \int_{0}^{t}{\delta _S(\tilde{x}(\tau ))d\tau }$ to
the projection of $\tilde{x}(t)$ on $S$.

The paper consists of two sections. In Section 2 we describe the analytic
method. It is similar to that of \cite{2}. The probabilistic approach is
presented in Section 3.

\head
{\S 2. The Analytic Method.}
\endhead

Let $\nu \in \Re ^d$ be a fixed vector with $\left\|\nu \right\|=1$ and $S$
be the hyperplane in $\Re ^d$ orthogonal to $\nu $
$$
S=\left\{x\in \Re ^d| (x,\nu )=0\right\}
$$

And let $x(t)$ be a continuous Markov process in $\Re ^d$ with transition
probability density $g_0(t,x,y)$
$$
g_0(t,x,y)=\left[det(2\pi tB)\right] ^{-1/2} exp\left\{-\frac{1}{2t}(B^{-1}
(y-x),y-x)\right\},  t>0, x\in \Re ^d, y\in \Re ^d
$$
where $B$ is the given positive definite symmetric operator in $\Re ^d$.

Let $q\in [-1,1],\alpha \in S$ be given. We are looking for a function
$u(t,x,\varphi)$ defined for $t>0, x\in \Re ^d, \varphi \in C_b(\Re ^d)$
such that

\roster

\item
$u(t,x,\varphi)$ satisfies the heat equation in the domain $t>0, x\not\in S$
$$
\frac{\partial  u}{\partial  t} = \frac{1}{2} Tr(BD^2u)
$$
If we fix a coordinates in $\Re ^d$ then $D^2u = (\frac{\partial  ^2u}
{\partial  x_i\partial  x_j})_{i,j=1}^d$ and this equation can be written
as follows
$$
\frac{\partial  u}{\partial  t} = \frac{1}{2} \sum_{i,j=1}^{d}{b_{ij}
\frac{\partial  ^2u}{\partial  x_i\partial  x_j}}
$$
where  $B=(b_{ij})_{i,j=1}^{d}$.

\item
$u(t,x,\varphi )$ satisfies the initial condition
$$
\lim_{t\downarrow 0}u(t,x,\varphi )=\varphi (x),  x\in \Re ^d, \varphi
\in C_b (\Re ^d)
$$

\item
For $x\in S$ and $t>0$ the relation
$$
u(t,x+,\varphi )=u(t,x-,\varphi )
$$
holds, where $u(t,x\pm ,\varphi )=\lim_{\varepsilon \downarrow 0}u(t,x\pm
\varepsilon \nu ,\varphi )$

\item
The relation
$$
\frac{1+q}{2}\frac{\partial  u(t,x+,\varphi )}{\partial  N}-\frac{1-q}{2}
\frac{\partial  u(t,x-,\varphi )}{\partial  N}+(\alpha ,\bigtriangledown u)=0,
$$
is valid for $x\in S$ and $t>0$ , where $\bigtriangledown u$ is gradient of
$u(t,x,\varphi )$, $N$ is the conormal vector to $S, N=B\nu $.

\endroster

Let us denote by $x_\nu (t)$ the projection of $x(t)$ on $\nu $, by $\pi _S$
the operator of orthogonal projection on $S$ such that $x^S(t)=\pi _Sx(t)$ is
the projection of $x(t)$ on $S$. If $\nu $ is not an eigen vector of the
operator $B$ then $x_\nu (t)$ and $x^S(t)$ are not independent processes and
we can not represent $g_0(t,x,y)$ as a product of $x_\nu (t)$ and $x^S(t)$
densities. But we may consider the following decomposition of $g_0(t,x,y)$
$$
g_0(t,x,y)=\frac{1}{\sqrt{2\pi t\sigma ^2}}\exp \left\{-\frac{(y_\nu -x_\nu )
^2}{2t\sigma ^2}\right\}
 \times
$$
$$
\times  \left[det(2\pi tB_S)\right] ^{-1/2} exp\left\{-\frac{1}{2t}(B^{-1}_S
(y_S-x_S-\frac{y_\nu -x_\nu }{\sigma ^2}b),y_S-x_S-\frac{y_\nu -x_\nu }{\sigma
^2}b)\right\}
$$
where $b=\pi _SB\nu , B_S=(\pi _SB^{-1}\pi _S)^{-1}, \sigma ^2=(B\nu ,\nu )$.

The first factor of this decomposition is the density of the process $x_\nu
(t)$, the second one is the conditional density of the process  $x^S(t)$ under
the condition $x_\nu (t)=y_\nu $. Denoting by $g^S(t,x,y)$ the second factor
we may write
$$
g_0(t,x,y)=\frac{1}{\sqrt{2\pi t\sigma ^2}}\exp \left\{-\frac{(y_\nu -x_\nu )^2}
{2t\sigma ^2}\right\}g^S(t,x,y),   t>0, x\in \Re ^d, y\in \Re ^d
$$

The problem $(1)-(4)$ may be solved in the following manner. Let us fix a
system of coordinates so that $\nu =e_d$ and make use of the Fourier
transformation with respect to $x_1,...,x_{d-1}$ and the Laplace transformation
with respect to $t$ to the function $u(t,x,\varphi )$. It will reduce this
problem to solving a pair of second order differential equations with constant
coefficients. Then free coefficients may be found from boundary conditions
$(3),(4)$. Returning to the original function we get the following assertion.

\proclaim{Theorem 1}
If $|q|\leq 1$ then there exists a continuous Markov process in $\Re ^d$ with
transition probability density $G(t,x,y)$ such that
$$
u(t,x,\varphi )=\int_{\Re ^d}{\varphi (y)G(t,x,y)dy},   t>0, x\in \Re ^d,
\varphi \in C_b(\Re ^d)
$$
where
$$
G(t,x,y)=\frac{1}{\sqrt{2\pi t\sigma ^2}}\left[\exp \left\{-\frac{(y_\nu -x_\nu )
^2}{2t\sigma ^2}\right\}-\exp \left\{-\frac{(|y_\nu |+|x_\nu |)^2}{2t\sigma ^2}
\right\}\right]g^S(t,x,y) +
$$
$$
+ \int_{0}^{\infty }{(q \sign y_\nu +1)\frac{\sigma ^2\theta +|x_\nu | +|y_\nu |}
{\sqrt{2\pi t^3\sigma ^2}}\exp \left\{-\frac{(\sigma ^2\theta +|x_\nu | +
|y_\nu |)^2}{2t\sigma ^2}\right\}g^S(t,x+\alpha \theta ,y)d\theta}
$$
\endproclaim

\remark{Remark} The condition $|q|\leq 1$ provides that $G(t,x,y)$ is positive.
\endremark

\head
{\S 3. The Probabilistic Method.}
\endhead

Let $w(t)$ be a given $d$-dimensional Wiener process, $B$ be a symmetric
positive definite operator in $\Re ^d$ and parameters $q\in [-1,1], \alpha
\in S$ be given. Consider the stochastic differential equation
$$
dx(t)=(qB\nu +\alpha )\delta _S(x(t))dt + B^{1/2}dw(t)\tag 1
$$
where $\delta _S(x)$ is the generalized function on $\Re ^d$ defined above.

When $\alpha =0$ then the solution of the equation (1) is known. It may be
called a multidimensional analogy to skew Brownian motion. Let us denote it
by $\tilde{x}(t)$. The density of $\tilde{x}(t)$ equals to
$$
\tilde{g}(t,x,y)=\frac{1}{\sqrt{2\pi t\sigma ^2}}\left[\exp \left\{-\frac
{y_\nu -x_\nu )^2}{2t\sigma ^2}\right\}+ q \sign y_\nu \exp \left\{-\frac{
(|y_\nu |+|x_\nu |)^2}{2t\sigma ^2}\right\}\right]g^S(t,x,y)
$$
where $t>0, x\in \Re ^d, y\in \Re ^d$

The equation (1) means that we have to find the process $x(t)$ that satisfies
the following properties. The projection of $x(t)$ on $\nu $ is skew Brownian
motion, the projection of $x(t)$ on $S$ is the Gaussian process in $S$ (with
mean $x^S+\frac{y_\nu -x_\nu }{\sigma ^2}b$ and correlation operator $tB_S$)
plus the process $\alpha \eta _t$, where $\eta _t$ is the functional of the
process $x(t)$ such that
$$
\eta _t=\int_{0}^{t}{\delta _S (x(\tau ))d\tau }
$$

\remark{Remark} Actually $\eta _t$ depends only on $x_\nu (t)$. Therefore we
may write $\eta _t=\int_{0}^{t}{\delta _S(\tilde{x}(\tau ))d\tau }$
\endremark

The functional $\eta _t$ is a nonnegative continuous homogeneous additive
functional, it increases at those instants of time for which the process
$\tilde{x}(t)$ hits the hyperplane $S$.

To solve the equation (1) we have to find the joint distribution of $\tilde{x}
(t)$ and $\eta _t$.

\proclaim{Lemma 1}
The joint distribution of $\tilde{x}(t)$ and $\eta _t$ has the form
$$
\P_x\left\{\tilde{x}(t)\in dy,\eta _t\in d\theta \right\}=\{\frac{\delta
(\theta )}{\sqrt{2\pi t\sigma ^2}}\left[\exp \left\{-\frac{(y_\nu -x_\nu )^2}
{2t\sigma ^2}\right\}-\exp \left\{-\frac{(|y_\nu |+|x_\nu |)^2}{2t\sigma ^2}
\right\}\right]+
$$
$$
+ (1+q \sign y_\nu )\frac{\sigma ^2\theta +|x_\nu |+|y_\nu |}{\sqrt{2\pi
t^3\sigma ^2}}\exp \left\{-\frac{(\sigma ^2\theta +|x_\nu | +|y_\nu |)^2}
{2t\sigma ^2}\right\}\}g^S(t,x,y)dyd\theta
$$
\endproclaim

\remark{Remark} $\P_x\left\{\tilde{x}(t)\in dy,\eta _t\in d\theta \right\}$
has the atom at the point $\theta =0$.
\endremark

\demo{Proof} Consider
$$
u(t,x,\lambda ,\mu )=\tilde{\E_x}e^{i(\tilde{x}(t),\mu )+i\lambda \eta _t},
t>0,x\in \Re ^d, \lambda \in \Re , \mu \in \Re ^d
$$
where $\tilde{\E_x}$ is the averaging operator corresponding $\tilde{g}(t,x,y)$.

The function $u(t,x,\lambda ,\mu )$ satisfies the following integral equation
(see \cite{3})
$$
u(t,x,\lambda ,\mu )=\int_{\Re ^d}{e^{i(y,\mu )}\tilde {g}(t,x,y)dy}+i\lambda
\int_{0}^{t}{d\tau \int_{S}{u(t-\tau ,z,\lambda ,\mu )\tilde{g}(\tau ,x,z)dz}}
\tag 2
$$

Note that the function $u(t,\cdot ,\lambda ,\mu )$ on the right hand side of
$(2)$ depends on its values on $S$. Then using Fourier-Laplace transformation
we can find $u(t,x,\lambda ,\mu )$ for $x\in S$
$$
u(t,x,\lambda ,\mu )=\int_{0}^{\infty }{d\theta \int_{\Re ^d}{e^{i(y,\mu )+
i\lambda \theta }(1+q \sign y_\nu )}}\times
$$
$$
\times \frac{1}{\sqrt{2\pi t\sigma ^2}}\exp \left\{-\frac{(\sigma ^2\theta
+|x_\nu | +|y_\nu |)^2}{2t\sigma ^2}\right\}g^S(t,x,y)dy,
$$
Substituting this relation in $(2)$ we will obtain
$$
u(t,x,\lambda ,\mu )=\int_{\Re ^d}\!\!{\frac{e^{i(y,\mu )}}{\sqrt{2\pi
t\sigma ^2}}\left[\exp \left\{\frac{-(y_\nu -x_\nu )^2}{2t\sigma ^2}\right\}
-\exp \left\{\frac{-(|y_\nu |+|x_\nu |)^2}{2t\sigma ^2}\right\}\right]g^S(t,x,y)
dy} +
$$
$$
+ \int_{0}^{\infty }{d\theta \int_{\Re ^d}{e^{i(y,\mu )+i\lambda \theta }}
(1+q \sign y_\nu )\frac{\sigma ^2\theta +|x_\nu |+|y_\nu |}{\sqrt{2\pi t^3
\sigma ^2}}}\times
$$
$$
\times \exp \left\{-\frac{(\sigma ^2\theta +|x_\nu | +|y_\nu |)^2}{2t\sigma
^2}\right\}g^S(t,x,y)dy
$$
where $t>0,x\in \Re ^d, \lambda \in \Re , \mu \in \Re ^d$. This gives us $
\P_x\left\{\tilde{x}(t)\in dy,\eta _t\in d\theta \right\}$.
\enddemo

If we know the joint distribution of $\tilde{x}(t)$ and $\eta _t$, then
the equation $(1)$ can be solved.

\proclaim{Theorem 2}
The solution of the equation $(1)$ is a continuous Markov process $x(t)$
with transition probability density $G(t,x,y)$.
\endproclaim

\demo{Proof}
We have
$$
\P_x\left\{x(t)\in dy \right\}=\int_{[0,+\infty )}{\P_x\left\{x(t)\in dy,
\eta _t\in d\theta \right\}}\tag 3
$$

When $\theta =0$ then $x(t)=\tilde{x}(t)$ and joint distribution $\P_x\left
\{x(t)\in dy,\eta _t\in d\theta \right\}$ is equal to $\P_x\left\{\tilde{x}
(t)\in dy,\eta _t\in d\theta \right\}$ at the point $\theta =0$. For each $
\theta \in (0,\infty )$ the process $x(t)$ is equal to $\tilde{x}(t)+\alpha
\theta $. Integrating $(3)$ on $\theta $ we obtain that
$$
\P_x\left\{x(t)\in dy \right\}=\frac{1}{\sqrt{2\pi t\sigma ^2}}\left[\exp
\left\{-\frac{(y_\nu -x_\nu )^2}{2t\sigma ^2}\right\}-\exp \left\{-\frac{
(|y_\nu |+|x_\nu |)^2}{2t\sigma ^2}\right\}\right]g^S(t,x,y)dy 
$$
$$
+ \int_{0}^{\infty }{(q \sign y_\nu +1)\frac{\sigma ^2\theta +|x_\nu |
+|y_\nu |}{\sqrt{2\pi t^3\sigma ^2}}\exp \left\{-\frac{(\sigma ^2\theta
+|x_\nu | +|y_\nu |)^2}{2t\sigma ^2}\right\}g^S(t,x+\alpha \theta ,y)
d\theta dy}
$$
It coincides with the result of the Theorem 1.
\enddemo

\enddocument